\documentclass[12pt]{article}
\usepackage{}
\usepackage{mathrsfs}
\usepackage[justification=centering]{caption}
\usepackage{booktabs}
\usepackage{multirow}
\usepackage{array}
\usepackage[numbers,sort&compress]{natbib}
\usepackage[colorlinks,
            linkcolor=blue,
            anchorcolor=green,
            citecolor=magenta
            ]{hyperref}
\usepackage{mathrsfs,amssymb,amsfonts}
\usepackage{graphicx}
\usepackage{amsmath,amssymb,latexsym,color}
\usepackage[mathscr]{eucal}
\usepackage[numbers]{natbib}
\usepackage{verbatim}
\newtheorem{thm}{Theorem}[section]
\newtheorem{definition}[thm]{Definition}

\newtheorem{lemma}[thm]{Lemma}

\newtheorem{corollary}[thm]{Corollary}
\newtheorem{example}{Example}
\newtheorem{remark}{Remark}[section]

\newcommand{\proof}{{\it Proof.\quad}}
\newcommand{\qed}{\hfill\Box\medskip}
\usepackage{CJK}
\usepackage{setspace}

\begin{document}

\title{\bf Hybrid fault diagnosis capability analysis of highly connected graphs
}

\author{
Yulong Wei\textsuperscript{a}\footnote{\footnotesize Corresponding author.\newline
\indent\indent{\em E-mail address:} weiyulong@tyut.edu.cn (Y. Wei).}
\quad
Rong-hua Li\textsuperscript{b}
\quad
Weihua Yang\textsuperscript{a}
\\
\\{\footnotesize
\textsuperscript{a}\em Department of Mathematics, Taiyuan University of Technology, Taiyuan, 030024, China}
\\{\footnotesize \textsuperscript{b}\em School of Computer Science {\rm \&} Technology, Beijing Institute of Technology, Beijing, 100081, China}}
\date{}
\maketitle

\setlength{\baselineskip}{24pt}

\noindent {\bf Abstract}\quad Zhu et al. [Theoret. Comput. Sci. 758 (2019) 1--8] introduced the $h$-edge tolerable diagnosability to measure the fault diagnosis capability of a multiprocessor system with faulty links. This kind of diagnosability is a generalization of the concept of traditional diagnosability. A graph is called a maximally connected graph if its minimum degree equals its vertex connectivity. It is well-known that many irregular networks are maximally connected graphs and the $h$-edge tolerable diagnosabilities of these networks are unknown, which is our motivation for research. In this paper, we obtain the lower bound of the $h$-edge tolerable diagnosability of a $t$-connected graph and establish the $h$-edge tolerable diagnosability of a maximally connected graph under the PMC model and the MM$^*$ model, which extends some results in [IEEE Trans. Comput. 23 (1974) 86--88], [IEEE Trans. Comput. 53 (2004) 1582--1590] and [Theoret. Comput. Sci. 796 (2019) 147--153].

\noindent {\bf Keywords}\quad Highly connected graph; Maximally connected graph; Fault diagnosability; PMC model; MM$^*$ model

\vskip0.6cm

\section{Introduction}
Processor failure has become an ineluctable event in a large-scale multiprocessor system. To keep the multiprocessor system performing its functions efficiently and economically, recognizing faulty processors correctly is a task of top priority. The process of recognizing faulty processors
in a multiprocessor system is called {\em fault diagnosis}, and the {\em diagnosability} of a system is the maximum number of faulty processors the system can recognize. The PMC model and the MM$^*$ model are two major models to investigate fault diagnosis in previous researches. The PMC model, proposed by Preparata, Metze and Chien \cite{Pre}, assumes that all adjacent processors of a system can test one another. The MM$^*$ model which is the development of the MM model \cite{Mae},
proposed by Sengupta and Dahbura \cite{Sen}, assumes that each processor has to test two processors if the processor is adjacent to the latter two processors. Some references related to fault diagnosis studies under the PMC model or MM$^*$ model can be seen in \cite{ChaC,Chen,Edd,Gu,Lin,Hua1,Lian,LinX,Wan,Wei,Wei186,Wei18,Wei188,Wei21,XuTH,XuT,Wei19,ZhaH,Zhan,ZhuL}.

In the real situation, both node and link faults can appear in a system. However, the traditional diagnosability for a multiprocessor system assumes that the system is without link faults. On the other hand, it is a natural question to ask how the diagnosability decreases if some links are missing for a multiprocessor system \cite{Wan}. To address the deficiency of the traditional diagnosability for a multiprocessor system and answer the above question, the concept of the $h$-edge tolerable diagnosability $t^e_h(G)$, introduced by Zhu et al. \cite{ZhuL}, generalizes the theories of diagnosability and can better measure the diagnosis capability of a multiprocessor system $G$. In fact, this diagnosability is the worst-case diagnosability when the number of faulty links of $G$ does not exceed $h$. Briefly, $t^e_h(G)$ is the minimum diagnosability of graphs $G-F_e$ which satisfy that $F_e\subseteq E(G)$ and $|F_e|\leq h$. Note that if a processor $u$ has no fault-free neighbors, it is impossible to determine whether $u$ is faulty or not in the fault diagnosis. Then $t^e_h(G)=0$ for $h\geq \delta(G)$, where $\delta(G)$ is the minimum degree of a graph $G$. Hence, a key issue for the $h$-edge tolerable diagnosability of a graph $G$ study is the case of $0\leq h\leq \delta(G)$.

In 2019, Zhu et al. \cite{ZhuL} determined the $h$-edge tolerable diagnosabilities of hypercubes under the PMC model and the MM$^*$ model. Wei and Xu \cite{Wei188,Wei19} established the $h$-edge tolerable diagnosabilities of  $k$-regular triangle-free graphs and balanced hypercubes under the PMC model and the MM$^*$ model. Recently, Lian et al. \cite{Lian} established the $h$-edge tolerable diagnosability of a $k$-regular $k$-connected graph under the PMC model and the MM$^*$ model. Zhang et al. \cite{Zhan} determined the $h$-edge tolerable diagnosabilities of triangle-free graphs under the PMC model and the MM$^*$ model, which extends the results of triangle-free regular graphs \cite{Wei19}. Zhang et al. \cite{ZhaH} determined the $h$-edge tolerable diagnosabilities of $k$-regular $2$-cn graphs under the PMC model for $h\leq k-5$.

A graph is called a maximally connected graph if its minimum degree equals its vertex connectivity. In this paper, we obtain the lower bound of the $h$-edge tolerable diagnosability of a $t$-connected graph and establish the $h$-edge tolerable diagnosability of a maximally connected graph under the PMC model and the MM$^*$ model, which provides a more precise characterization for the fault diagnosis capability of networks and generalizes some results in \cite{ChaC,Hak,Lian}.

The remainder of this paper is organized as follows. Some terminology and preliminaries are introduced in Section \ref{2}. The main results are given in Sections \ref{3} and \ref{4}.  Finally, we concludes the paper in Section \ref{6}.

\section{Terminology and preliminaries}\label{2}
A {\em graph} $G =\big(V(G), E(G)\big)$ is used to represent a system (or a network), where each vertex of $G$ represents a processor and each edge of $G$ represents a link. The {\em connectivity} $\kappa(G)$ is the minimum cardinality of all vertex subsets $S\subseteq V(G)$ satisfying that $G-S$ is disconnected or trivial. A graph $G$ is said to be {\em $t$-connected}, if $\kappa(G)\geq t$. The {\em neighborhood} $N_G(v)$ of a vertex $v$ in $G$ is the set of vertices adjacent to $v$. We call $\min_{v\in V(G)}\{|N_G(v)|\}$ the {\em minimum degree} of a graph $G$, denoted by $\delta(G)$. A graph $G$ is said to be {\em $t$-regular} (or {\em regular}), if $|N_G(v)|=t$ for any vertex $v$ of $G$. We refer readers to \cite{Bon} for terminology and notation unless stated otherwise.

The concept of the traditional diagnosability of a graph is presented as follows.
\begin{definition}[\cite{Dah}]\label{D3}
A graph $G=(V, E)$ of $n$ vertices is $t$-diagnosable if all faulty vertices can be detected without replacement, provided that the number of faults does not exceed $t$. The diagnosability $t(G)$ of a graph $G$ is the maximum value of $t$ such that $G$ is $t$-diagnosable.
\end{definition}

\begin{figure}[hptb]
 \centering
  \includegraphics[width=10cm]{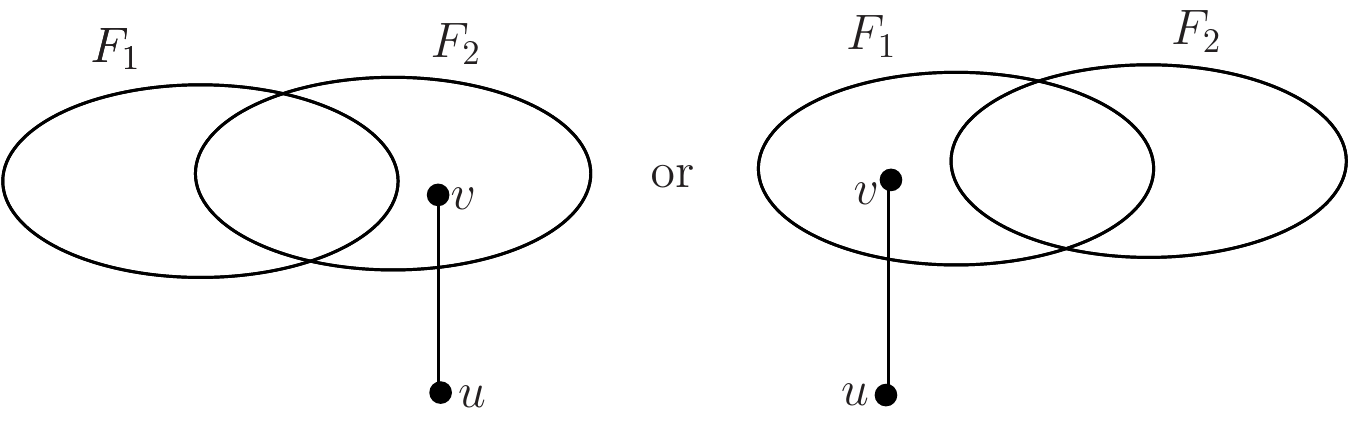}\\
  \caption{ The illustration of Lemma \ref{L01}.
}\label{L11}
\end{figure}

\begin{figure}[hptb]
  \centering
  \includegraphics[width=6cm]{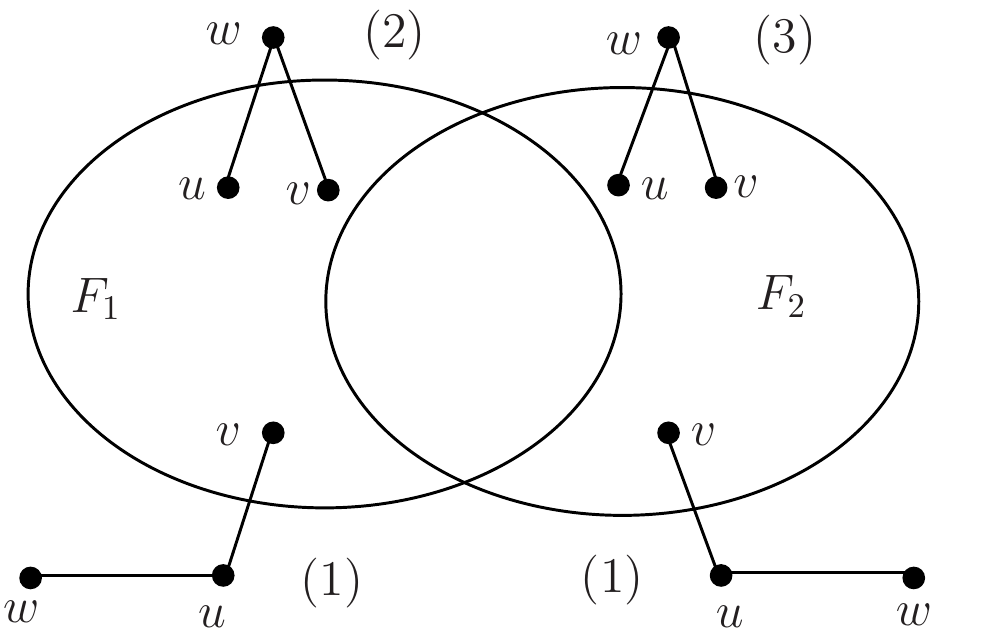}\\
  \caption{ The illustration of Lemma \ref{L02}.
}\label{L22}
\end{figure}

For any two sets $A$ and $B$, we use $A-B$ to denote a set obtained by removing all elements of $B$ from $A$. The {\em symmetric difference} of two sets $F_1$ and $F_2$ is defined as the set $F_1\bigtriangleup F_2$ $=(F_1-F_2)\cup (F_2-F_1)$. The following lemmas give necessary and sufficient conditions for a graph to be $t$-diagnosable under the PMC model and the MM$^*$ model.
\begin{lemma}[\cite{Dah}]\label{L01}
A graph $G =(V, E)$ is $t$-diagnosable under the PMC model if and only if for any two distinct subsets $F_1$ and $F_2$ of $V$ with $|F_1|\leq t$ and $|F_2|\leq t$, there exists an edge from $V-(F_1\cup F_2)$ to $F_1\bigtriangleup F_2$ {\rm(}see Figure \ref{L11} {\rm)}.
\end{lemma}
\begin{lemma}[\cite{Sen}]\label{L02}
A graph $G =(V, E)$ is $t$-diagnosable under the MM$^*$ model if and only if for any two distinct subsets $F_1$ and $F_2$ of $V$ with $|F_1|\leq t$ and $|F_2|\leq t$, at least one of the following conditions is satisfied {\rm(}see Figure \ref{L22} {\rm)}:
\begin{enumerate}
\item[{\rm (1)}] There are two vertices $u, w\in V-(F_1\cup F_2)$ and there is a vertex $v\in F_1\bigtriangleup F_2$
such that $uv\in E$ and $uw\in E$.

\item[{\rm (2)}] There are two vertices $u, v\in F_1-F_2$ and there is a vertex $w\in V-(F_1\cup F_2)$
such that $uw\in E$ and $vw\in E$.

\item[{\rm (3)}] There are two vertices $u, v\in F_2-F_1$ and there is a vertex $w\in V-(F_1\cup F_2)$
such that $uw\in E$ and $vw\in E$.
\end{enumerate}
\end{lemma}

We call sets $F_1$ and $F_2$ {\em distinguishable} under the PMC (resp. MM$^*$) model if they satisfy the condition of Lemma \ref{L01} (resp. at least one of the conditions of Lemma \ref{L02}). Otherwise, $F_1$ and $F_2$ are said to be {\em indistinguishable}.

To better adapt to the real circumstances that link faults may happen \cite{ZhuL}, Zhu et al. introduced the definition of the $h$-edge tolerable diagnosability of graphs as follows.

\begin{definition}\label{D1}
Given a diagnosis model and a graph $G$, $G$ is $h$-edge tolerable $t$-diagnosable under the
diagnosis model if for any edge subset $F_e$ of $G$ with $|F_e|\leq h$, the graph $G-F_e$ is $t$-diagnosable under the
diagnosis model. The $h$-edge tolerable diagnosability of $G$, denoted as $t_h^e(G)$, is the maximum integer $t$ such that $G$ is $h$-edge tolerable $t$-diagnosable.
\end{definition}
Obviously, for a graph $G$, the $h$-edge tolerable diagnosability is the traditional diagnosability when $h=0$. We have $t_0^e(G)=t(G)$.

A family of paths in $G$ is said to be {\em internally-disjoint} if no vertex of $G$ is an internal vertex of more than one path of the family. The following lemmas are important to the proof of our main results.

\begin{lemma}[Whitney(1932) \cite{Bon}]\label{LL02}
A graph $G$ with at least three vertices is $2$-connected if and only if any two vertices of $G$ are connected by at least two internally-disjoint paths.
\end{lemma}

\begin{lemma}[\cite{Lian}]\label{LL0}
Let $G=(V, E)$ be a connected graph and $S\subseteq E$. If $|S|\leq \kappa(G)$, then $\kappa(G-S)\geq \kappa(G)-|S|$.
\end{lemma}

\begin{lemma}[\cite{Zhan}]\label{LL1}
Let $G=(V, E)$ be a connected graph with minimum degree $\delta(G)$. Then $t_h^e(G)\leq \delta(G)-h$ under the PMC model and the MM$^*$ model for $0\leq h \leq \delta(G)$. Particularly, $t_{\delta(G)}^e(G)=0$ under both the PMC model and the MM$^*$ model.
\end{lemma}

\section{Hybrid fault diagnosis capability analysis of highly connected graphs under the PMC model}\label{3}
In this section, we will discuss the $h$-edge tolerable
diagnosability of a highly connected graph under the PMC model.

First, we give a lower bound of the $h$-edge tolerable diagnosability of a $t$-connected graph under the PMC model.

\begin{thm}\label{LL2}
Let $G=(V, E)$ be a $t$-connected graph with $|V|\geq 2(t-h)+1$. Then $t_h^e(G)\geq t-h$
under the PMC model for $0\leq h \leq t$.
\end{thm}
\proof
If $h=t$, then $t_h^e(G)\geq 0=t-h$ holds obviously.

Now, we assume that $h \leq t-1$. For an arbitrary edge subset $F_e\subseteq E$ with $|F_e|\leq h$, suppose that there exist two distinct vertex subsets $F_1, F_2\subseteq V$ such that $F_1$ and $F_2$ are indistinguishable in $G-F_e$ under the PMC model. We will prove this theorem by showing that $|F_1|\geq t-h+1$ or $|F_2|\geq t-h+1$ for $0\leq h \leq t-1$.

If $|F_1\cap F_2|\geq t-h$, then $|F_1|\geq t-h+1$ or $|F_2|\geq t-h+1$.

Suppose that $|F_1\cap F_2|\leq t-h-1$. Then by Lemma \ref{LL0}, $$\kappa(G-F_e)\geq \kappa(G)-|F_e|\geq t-h>t-h-1\geq|F_1\cap F_2|. $$ Therefore, $G-F_e-(F_1\cap F_2)$ is connected.

If $V=F_1\cup F_2$, then $|F_1\cup F_2|=|V|\geq 2(t-h)+1$. Thus, $|F_1|\geq t-h+1$ or $|F_2|\geq t-h+1$. Otherwise, $V-(F_1\cup F_2)\neq\emptyset$. Since $G-F_e-(F_1\cap F_2)$ is connected, there is an edge between $F_1\bigtriangleup F_2$ and $V-(F_1\cup F_2)$ in $G-F_e$, a contradiction by Lemma \ref{L01}. $\qed$

Let $h=0$. By Theorem \ref{LL2}, we can obtain the following result.
\begin{corollary}[\cite{Hak}]\label{cor4}
Let $G=(V, E)$ be a $t$-connected network with $N$ nodes and $t\geq2$. $G$ is $t$-diagnosable
under the PMC model if $N\geq 2t+1$.
\end{corollary}

Note that a maximally connected graph $G$ is $\delta(G)$-connected. By Lemma \ref{LL1} and Theorem \ref{LL2}, we obtain the following result.
\begin{thm}\label{main1}
Let $G=(V, E)$ be a maximally connected graph with $|V|\geq 2(\delta(G)-h)+1$. Then $t_h^e(G)=\delta(G)-h$
under the PMC model for $0\leq h \leq \delta(G)$.
\end{thm}

Note that a $k$-regular $k$-connected graph is a maximally connected graph. By Theorem \ref{main1}, we immediately obtain the following result.
\begin{corollary}[\cite{Lian}]\label{cor1}
Let $G=(V, E)$ be a $k$-regular $k$-connected graph with $|V|\geq 2(k-h)+1$. Then $t_h^e(G)=k-h$
under the PMC model for $0\leq h \leq k$.
\end{corollary}

\section{Hybrid fault diagnosis capability analysis of highly connected graphs under the MM$^*$ model}\label{4}
In this section, we will discuss the $h$-edge tolerable
diagnosability of a highly connected graph under the MM$^*$ model.

\begin{figure}[hptb]
  \centering
  \includegraphics[width=12cm]{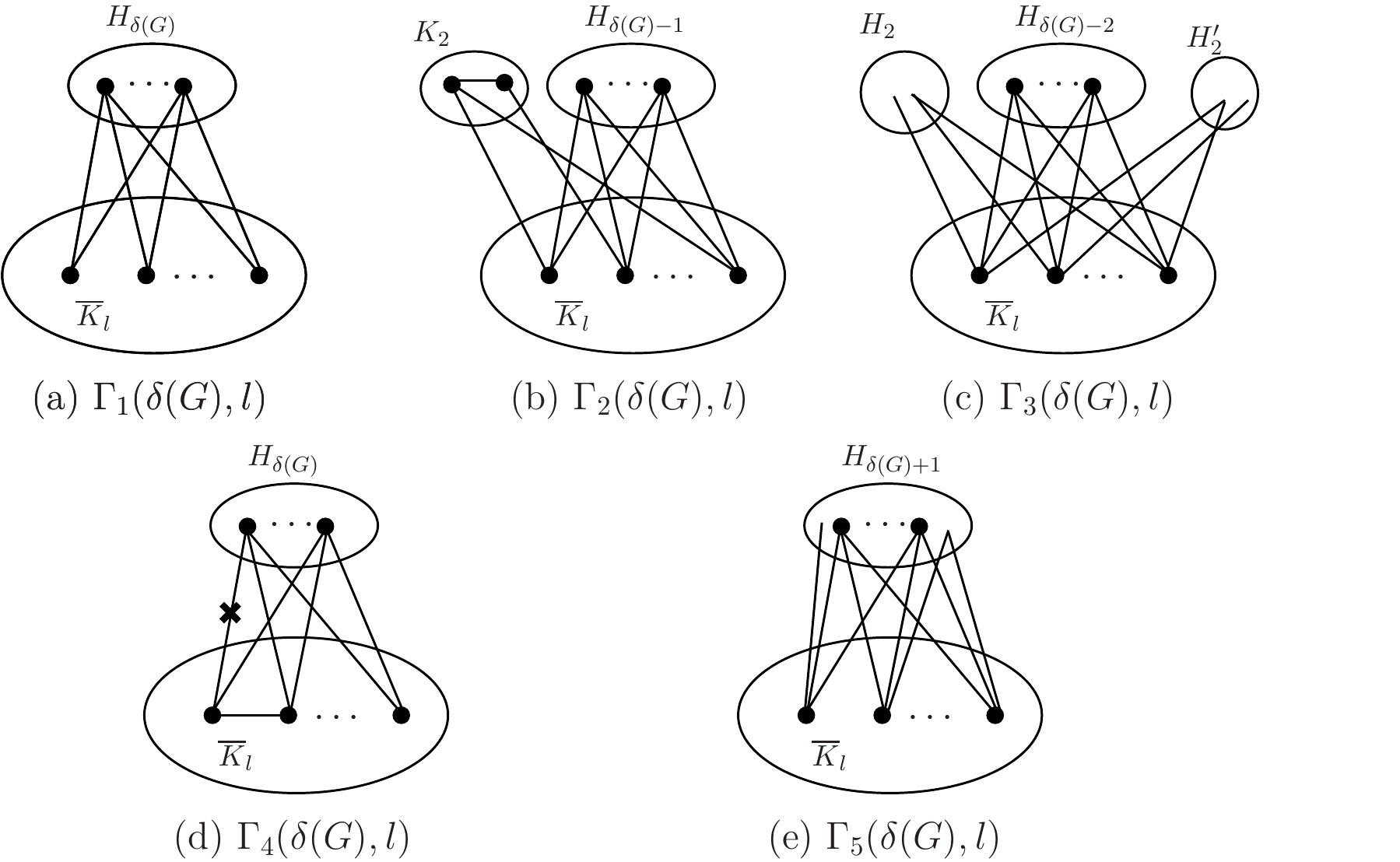}\\
  \caption{The graphs $\Gamma_i(\delta(G), l)$ for $1\leq i\leq5$. }\label{11}
\end{figure}

In the following statements, we use $\overline{G}$ to denote the {\em complement} graph of a simple graph $G$, whose vertex set is $V(G)$ and whose edges are the pairs of nonadjacent vertices of $G$. By starting with a disjoint union of two graphs $G$ and $H$ (i.e., $G\cup H$) and adding edges joining every vertex of $G$ to every vertex of $H$, one obtains
the {\em join } of $G$ and $H$, denoted by $G\vee H$ \cite{Bon}.
For two graphs $G$ and $H$, we use $G\ast_{r} H$ to denote a graph with the vertex set $V(G)\cup V(H)$ and the edge set $E(G)\cup E(H)\cup E[V(G), V(H)]$, where $E[V(G), V(H)]$ is an edge subset and any edge in $E[V(G), V(H)]$ has one endpoint in $V(G)$ and the other endpoint in $V(H)$. Clearly, $0\leq|E[V(G), V(H)]|\leq|V(G)|\cdot|V(H)|$. We also use $G\ast_{1} H$ to denote a graph with the vertex set $V(G)\cup V(H)$ and the edge set $E(G)\cup E(H)\cup E[V(G), V(H)]$ such that any vertex of $G$ is adjacent to only one vertex of $H$. A spanning subgraph of a complete graph $K_{m}$ is denoted by $H_{m}$ (or $H_{m}'$). Given a graph $G$ with $\delta(G)\geq3$, let a collection of graphs $$\mathcal{F}(\delta(G))=\{\Gamma_i(\delta(G), l)\mid l\geq \delta(G)+1, 1\leq i\leq4\}\cup\{\Gamma_5(\delta(G), l)\mid l\geq \delta(G)+2\},$$ where $\Gamma_1(\delta(G), l)=H_{\delta(G)}\vee \overline{K}_l$,  $\Gamma_2(\delta(G), l)=(H_{\delta(G)-1}\ast_{r}K_2)\cup(H_{\delta(G)-1}\vee \overline{K}_l)\cup(\overline{K}_l\ast_{1}K_2)$, $\Gamma_3(\delta(G), l)=(\overline{K}_l\ast_{1}H_2)\cup(H_{\delta(G)-2}\vee \overline{K}_l)\cup(\overline{K}_l\ast_{1}H_2')\cup(H_2\ast_{r}H_{\delta(G)-2}\ast_{r}H_2')$, $\Gamma_4(\delta(G), l)=(H_{\delta(G)}\vee \overline{K}_l)+e-E_0$ with $e=u_1u_2$, $u_1, u_2\in V(\overline{K}_l)$, $E_0\subseteq E[V(H_{\delta(G)}), \{u_1, u_2\}]$ and $|E_0\cap E[V(H_{\delta(G)}), \{u_i\}]|\leq1$ for $1\leq i\leq2$, and $\Gamma_5(\delta(G), l)=(V(H_{\delta(G)+1})\cup V(\overline{K}_l), E(H_{\delta(G)+1})\cup\{uv\mid u\in V(\overline{K}_l), v\in V(H_{\delta(G)+1}), \delta(G)\leq|E[\{u\}, V(H_{\delta(G)+1})]|\leq\delta(G)+1\})$ (see Figure \ref{11}).

Let $C(G)$ be the maximum number of common neighbors of any two vertices in the graph $G$. We give some properties of graphs in $\mathcal{F}(\delta(G))$ as follows.

\begin{lemma}\label{Pro}
Given a graph $G$ with $\delta(G)\geq3$, the graphs in $\mathcal{F}(\delta(G))$ satisfy the following properties:
\begin{description}
  \item[(1)] The graphs in $\mathcal{F}(\delta(G))$ are all irregular;
  \item[(2)] $\min\{C(F)\mid F\in \mathcal{F}(\delta(G))\}\geq \delta(G)-1$.
\end{description}
\end{lemma}
\proof (1) Suppose $x\in V(H_{\delta(G)-i+1})$ and $y\in V(\overline{K}_l)$ for $1\leq i\leq3$. Since $|N_{\Gamma_i(\delta(G), l)}(x)|\geq l\geq \delta(G)+1>\delta(G)=|N_{\Gamma_i(\delta(G), l)}(y)|$, we know the graph $\Gamma_i(\delta(G), l)$ is irregular for $1\leq i\leq3$.

Since $\delta(G)\geq3>2\geq |E_0|$, there exists a vertex $x\in V(H_{\delta(G)})$ such that $V(\overline{K}_l)\subseteq N_{\Gamma_4(\delta(G), l)}(x)$. Pick a vertex $y\in V(\overline{K}_l)-\{u,v\}$. Note that $|N_{\Gamma_4(\delta(G), l)}(x)|\geq l\geq \delta(G)+1>\delta(G)=|N_{\Gamma_4(\delta(G), l)}(y)|$. We know that the graph $\Gamma_4(\delta(G), l)$ is irregular.

Assume to the contrary that $\Gamma_5(\delta(G), l)$ is a regular graph. Note that there exists a vertex
$x\in V(H_{\delta(G)+1})$ such that
$$|N_{\Gamma_5(\delta(G), l)}(x)|\geq \dfrac{|E[V(H_{\delta(G)+1}),V(\overline{K}_l)]|}{|V(H_{\delta(G)+1})|}\geq\dfrac{\delta(G)|V(\overline{K}_l)|}{|V(H_{\delta(G)+1})|}\geq \dfrac{\delta(G)(\delta(G)+2)}{\delta(G)+1}>\delta(G). $$
If $|N_{\Gamma_5(\delta(G), l)}(x)|=\delta(G)+1$, then $|N_{\Gamma_5(\delta(G), l)}(y)|=\delta(G)+1$ for any $y\in V(\overline{K}_l)$. Thus there exists a vertex
$z\in V(H_{\delta(G)+1})$ such that
$$|N_{\Gamma_5(\delta(G), l)}(z)|\geq \dfrac{|E[V(H_{\delta(G)+1}),V(\overline{K}_l)]|}{|V(H_{\delta(G)+1})|}= \dfrac{(\delta(G)+1)|V(\overline{K}_l)|}{\delta(G)+1}=l\geq\delta(G)+2,  $$
which contradicts that $\Gamma_5(\delta(G), l)$ is a regular graph. If $|N_{\Gamma_5(\delta(G), l)}(x)|\geq\delta(G)+2$, then $|N_{\Gamma_5(\delta(G), l)}(y)|\geq\delta(G)+2$ for any $y\in V(\overline{K}_l)$, which contradicts $N_{\Gamma_5(\delta(G), l)}(y)\subseteq V(H_{\delta(G)+1})$.

(2) Suppose $x, y\in V(H_{\delta(G)})$ (see Figure \ref{11}(a)). Since $|N_{\Gamma_1(\delta(G), l)}(x)\cap N_{\Gamma_1(\delta(G), l)}(y)|\geq|V(\overline{K}_l)|\geq \delta(G)+1$, we have $C(\Gamma_1(\delta(G), l))\geq \delta(G)+1$.

Note that $|N_{\Gamma_i(\delta(G), l)}(z)|=\delta(G)$ for any $z\in V(\overline{K}_l)$ and $i\in\{2,3\}$ (see Figure \ref{11}(b) and Figure \ref{11}(c)). Since $l\geq \delta(G)+1>2$, there exist two distinct vertices $x, y\in V(\overline{K}_l)$ such that $|N_{\Gamma_i(\delta(G), l)}(x)\cap N_{\Gamma_i(\delta(G), l)}(y)|\geq\delta(G)-i+2$. Thus, we have $C(\Gamma_i(\delta(G), l))\geq \delta(G)-i+2$, where $i\in\{2,3\}$.

Since $|V(H_{\delta(G)})|=\delta(G)\geq3>2\geq |E_0|$, there exist $x, y\in V(H_{\delta(G)})$ such that $V(\overline{K}_l)\subseteq N_{\Gamma_4(\delta(G), l)}(x)$ and $|N_{\Gamma_4(\delta(G), l)}(y)\cap V(\overline{K}_l)|\geq l-1\geq\delta(G)$ (see Figure \ref{11}(d)). Then $C(\Gamma_4(\delta(G), l))\geq \delta(G)$.

Note that $|N_{\Gamma_5(\delta(G), l)}(z)|\geq\delta(G)$ and $N_{\Gamma_5(\delta(G), l)}(z)\subseteq V(H_{\delta(G)+1})$ for any $z\in V(\overline{K}_l)$ (see Figure \ref{11}(e)). Since $l\geq \delta(G)+2>|V(H_{\delta(G)+1})|$, there exist two distinct vertices $x, y\in V(\overline{K}_l)$ such that $|N_{\Gamma_5(\delta(G), l)}(x)\cap N_{\Gamma_5(\delta(G), l)}(y)|\geq\delta(G)$. Thus, we have $C(\Gamma_5(\delta(G), l))\geq \delta(G)$.

As mentioned above, we get the desired results. $\qed$

\begin{remark}\label{r1}
According to the proof of the Lemma \ref{Pro}, if $\delta(G)\geq4$, then $\min\{C(F)\mid F\in \mathcal{F}(\delta(G))\}\geq\delta(G)$.
\end{remark}

Now, we give a lower bound of the $h$-edge tolerable diagnosability of a $t$-connected graph $G$ under the MM$^*$ model. For a given vertex $x\in V(G)$, we use $E(x)$ to denote the edges incident with $x$ in $G$.

\begin{thm}\label{LL3}
Let $G=(V, E)$ be a $t$-connected graph with $|V|\geq 2(t-h)+3$ and $t\geq3$. If $G\notin \mathcal{F}(\delta(G))$, then $t_h^e(G)\geq t-h$ under the MM$^*$ model for $0\leq h \leq \lfloor\dfrac{t-1}{2}\rfloor$.
\end{thm}
\proof
For an arbitrary edge subset $F_e\subseteq E$ with $|F_e|\leq h$, suppose that there exist two distinct vertex subsets $F_1, F_2\subseteq V$ such that $F_1$ and $F_2$ are indistinguishable in $G-F_e$ under the MM$^*$ model. We will prove the lemma by showing that $|F_1|\geq t-h+1$ or $|F_2|\geq t-h+1$ for $0\leq h \leq\lfloor\dfrac{t-1}{2}\rfloor$ and $t\geq3$.

If $|F_1\cap F_2|\geq t-h$, then $|F_1|\geq t-h+1$ or $|F_2|\geq t-h+1$.

Suppose that $|F_1\cap F_2|\leq t-h-1$. Then by Lemma \ref{LL0}, $$\kappa(G-F_e)\geq \kappa(G)-|F_e|\geq t-h>t-h-1\geq|F_1\cap F_2|. $$ Therefore, $G-F_e-(F_1\cap F_2)$ is connected.

If $V=F_1\cup F_2$, then $|F_1\cup F_2|=|V|\geq 2(t-h)+3$. Thus, $|F_1|\geq t-h+1$ or $|F_2|\geq t-h+1$. Now, we assume that $V-(F_1\cup F_2)\neq\emptyset$.

\textbf{Claim 1}: $V-(F_1\cup F_2)$ is an independent set of $G-F_e$.

Otherwise, $E(G[V-(F_1\cup F_2)])\neq\emptyset$. Since $G-F_e-(F_1\cap F_2)$ is connected, there exist three distinct vertices $x, y\in V-(F_1\cup F_2)$ and $w\in F_1\bigtriangleup F_2$ such that $xy, yw\in E(G-F_e)$. By Lemma \ref{L02}, $F_1$ and $F_2$ are distinguishable in $G-F_e$ under the MM$^*$ model, which is a contradiction.

\begin{figure}[hptb]
  \centering
  \includegraphics[width=8cm]{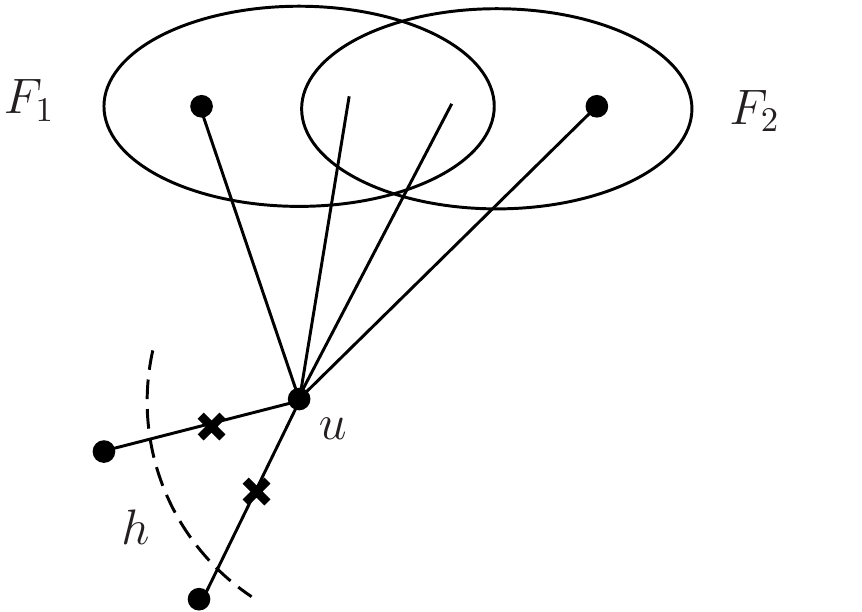}\\
  \caption{Illustration of $|F_1\cap F_2|\geq t-h-2$. }\label{22}
\end{figure}

Pick a vertex $u\in V-(F_1\cup F_2)$. Note that $F_1$ and $F_2$ are indistinguishable in $G-F_e$ under the MM$^*$ model.
By Claim 1, $|F_1\cap F_2|\geq|N_{G-F_e}(u)|-|N_{G-F_e}(u)\cap(F_1-F_2)|-|N_{G-F_e}(u)\cap(F_2-F_1)|\geq |N_G(u)|-h-2\geq\delta(G)-h-2\geq t-h-2$ (see Figure \ref{22}).

\textbf{Case 1}: $|F_1\cap F_2|=t-h-2$.

In this case, we have $|F_e|=h$, $F_e\subseteq E(u)$ and $|N_G(u)|=\delta(G)=t$. Thus, $|N_{G-F_e}(u)\cap(F_1-F_2)|=|N_{G-F_e}(u)\cap(F_2-F_1)|=1$.
If $|F_1-F_2|\geq3$ or $|F_2-F_1|\geq3$, then $|F_1|\geq t-h+1$ or $|F_2|\geq t-h+1$.
Now, we assume that $1\leq|F_1-F_2|\leq2$ and $1\leq|F_2-F_1|\leq2$.

\textbf{Claim 2}: $|V-(F_1\cup F_2)|\geq3$.

Note that $|V-(F_1\cup F_2)|=|V|-|F_1\cap F_2|-|F_1\bigtriangleup F_2|\geq 2(t-h)+3-(t-h-2)-4=t-h+1\geq t-\lfloor\dfrac{t-1}{2}\rfloor+1\geq3$. Thus, Claim 2 holds.

By Claim 2, suppose $v\in V-(F_1\cup F_2)-\{u\}$. Since $F_1$ and $F_2$ are indistinguishable in $G-F_e$ under the MM$^*$ model, we have $|N_{G-F_e}(v)\cap(F_1-F_2)|\leq1$ and $|N_{G-F_e}(v)\cap(F_2-F_1)|\leq1$.

If $E(u)\cap E(v)=\emptyset$ and $h\geq1$, then $|N_{G-F_e}(v)\cap(F_1\cap F_2)|\geq|N_{G-F_e}(v)|-2=|N_G(v)|-2\geq\delta(G)-2$,
which contradicts $|N_{G-F_e}(v)\cap(F_1\cap F_2)|\leq|F_1\cap F_2|=\delta(G)-h-2\leq\delta(G)-3$.

If $E(u)\cap E(v)\neq\emptyset$ and $h\geq2$, then $|N_{G-F_e}(v)\cap(F_1\cap F_2)|\geq|N_{G-F_e}(v)|-2\geq(|N_G(v)|-1)-2\geq\delta(G)-3$,
which contradicts $|N_{G-F_e}(v)\cap(F_1\cap F_2)|\leq|F_1\cap F_2|=\delta(G)-h-2\leq\delta(G)-4$.

If $E(u)\cap E(v)\neq\emptyset$ and $h=1$, then by Claim 2, there exists a vertex $w\in V-(F_1\cup F_2)-\{u,v\}$ such that $E(w)\cap F_e=\emptyset$. Therefore, $\delta(G)-1-2=|F_1\cap F_2|\geq|N_{G-F_e}(w)\cap(F_1\cap F_2)|=|N_{G}(w)\cap(F_1\cap F_2)|\geq\delta(G)-2$, a contradiction.

Now, we consider the case of $h=0$. Then $|F_1\cap F_2|=\delta(G)-2$ and $|V-(F_1\cup F_2)|=|V|-|F_1\cap F_2|-|F_1\bigtriangleup F_2|\geq 2\delta(G)+3-(\delta(G)-2)-4=\delta(G)+1$.
Note that $\delta(G)\leq|F_1\cup F_2|\leq \delta(G)+2$. Thus, we distinguish the following cases.

\textbf{Case 1.1}: $|F_1\cup F_2|=\delta(G)$.

By Claim 1, we have $G$ is isomorphic to $\Gamma_1(\delta(G), l)$ for some graph $H_{\delta(G)}$ and some integer $l$, a contradiction.

\textbf{Case 1.2}: $|F_1\cup F_2|=\delta(G)+1$.

Without loss of generality, we assume that $F_1-F_2=\{w_1, w_2\}$ and $F_2-F_1=\{w\}$. Since $G$ is a $t$-connected graph, we have that $\kappa(G)-|F_1\cap F_2|\geq t-(t-2)=2$. Then $G-(F_1\cap F_2)$ is $2$-connected. If $w_1w_2\notin E(G)$, then by Lemma \ref{LL02}, there exist two internally-disjoint paths $P$ and $Q$ from $w_1$ to $w_2$ in $G-(F_1\cap F_2)$.
\begin{figure}[hptb]
  \centering
  \includegraphics[width=8cm]{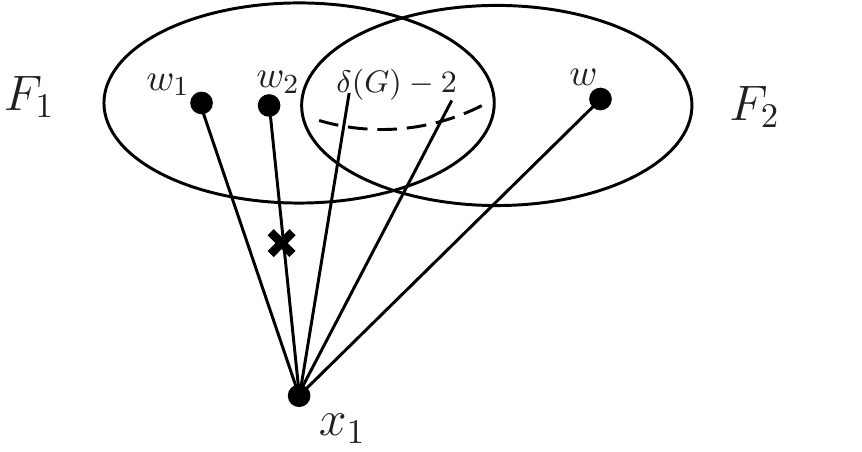}\\
  \caption{Illustration of $w_1w_2\in E(G)$. }\label{44}
\end{figure}
Without loss of generality, we assume that $w\notin V(P)$. Suppose that $P=(w_1,x_1,\ldots,x_m,w_2)$. Then $x_i\in V-(F_1\cup F_2)$ for $1\leq i\leq m$. If $m\geq2$, then $x_{m-1}x_m\in E(G)$, which contradicts Claim 1. Otherwise, $m=1$ and $x_1w_i\in E(G)$ for $1\leq i\leq2$, which contradicts that $F_1$ and $F_2$ are indistinguishable under the MM$^*$ model (see Figure \ref{44}). Thus, we know that $w_1$ is adjacent to $w_2$ and $G$ is isomorphic to $\Gamma_2(\delta(G), l)$ for some graph $H_{\delta(G)-1}$ and some integer $l$, which contradicts $G\notin \mathcal{F}(\delta(G))$.

\textbf{Case 1.3}: $|F_1\cup F_2|=\delta(G)+2$.

\begin{figure}[hptb]
  \centering
  \includegraphics[width=8cm]{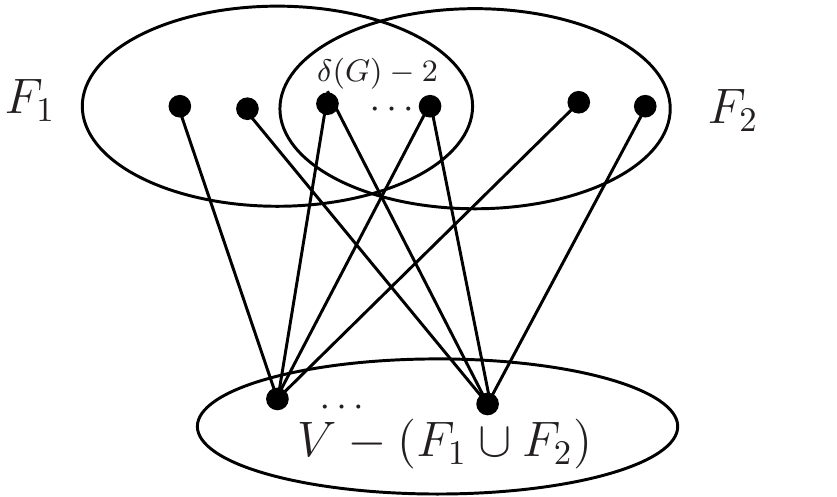}\\
  \caption{Illustration of the proof of Case 1.3.}\label{55}
\end{figure}

In this case, $|F_1-F_2|=|F_2-F_1|=2$ (see Figure \ref{55}). Note that $F_1$ and $F_2$ are indistinguishable under the MM$^*$ model. By Claim 1, we have $G$ is isomorphic to $\Gamma_3(\delta(G), l)$ for some graphs $H_{\delta(G)-2}$, $H_{2}$, $H'_{2}$ and some integer $l$, a contradiction.

\textbf{Case 2}: $|F_1\cap F_2|=t-h-1$.

Since $F_1\neq F_2$, without loss of generality, we assume that $F_2-F_1\neq\emptyset$. If $|F_1-F_2|\geq2$ or $|F_2-F_1|\geq2$, then $|F_1|\geq t-h+1$ or $|F_2|\geq t-h+1$.
Now, assume that $|F_1-F_2|\leq1$ and $|F_2-F_1|=1$. Thus, $|F_1\cup F_2|=|F_1\cap F_2|+|F_1\bigtriangleup F_2|\leq (t-h-1)+2=t-h+1$ and so $|V-(F_1\cup F_2)|\geq 2(t-h)+3-(t-h+1)=t-h+2\geq4$.

\textbf{Case 2.1}: $F_1-F_2=\emptyset$.

Note that $N_{G-F_e}(w)\subseteq F_1\cup F_2$ for any $w\in V-(F_1\cup F_2)$. Then, $\delta(G)-h\geq t-h=|F_1\cup F_2|\geq|N_{G-F_e}(w)|=|E(w)-F_e|\geq|E(w)|-|F_e|\geq\delta(G)-h$. Therefore, $t=\delta(G)$ and $F_e\subseteq E(w)$ for any $w\in V-(F_1\cup F_2)$. Thus, $F_e\subseteq \cap_{w\in V-(F_1\cup F_2)} E(w)$. Since $|V-(F_1\cup F_2)|\geq3$, we have $F_e=\emptyset$. Hence, $|F_1\cup F_2|=\delta(G)$ and $|V-(F_1\cup F_2)|\geq 2\delta(G)+3-\delta(G)=\delta(G)+3$. By Claim 1, we know that $G$ is isomorphic to $\Gamma_1(\delta(G), l)$ for some graph $H_{\delta(G)}$ and some integer $l$, which is a contradiction.

\textbf{Case 2.2}: $F_1-F_2\neq\emptyset$.

In this case, $|F_1-F_2|=|F_2-F_1|=1$ and $|F_1\cup F_2|=t-h+1$. If there exists a vertex $x\in V-(F_1\cup F_2)$ such that $|E(x)\cap F_e|\leq h-2$, then $|N_{G-F_e}(x)|=|E(x)|-|E(x)\cap F_e|\geq\delta(G)-(h-2)>\delta(G)-h+1\geq t-h+1=|F_1\cup F_2|$, which contradicts $N_{G-F_e}(x)\subseteq F_1\cup F_2$. Now, we assume that $|E(x)\cap F_e|\geq h-1$ for any vertex $x\in V-(F_1\cup F_2)$.

Note that $|V-(F_1\cup F_2)|\geq2(t-h)+3-(t-h+1)=t-h+2\geq (2h+1)-h+2=h+3$. If there exists a vertex $w\in V-(F_1\cup F_2)$ such that $|E(w)\cap F_e|=h$, then suppose $ww_i\in E(w)\cap F_e$ for $1\leq i\leq h$. Since $|V-(F_1\cup F_2)|\geq h+3>h+1=|\{w,w_1,\ldots,w_h\}|$, there exists a vertex $z\in V-(F_1\cup F_2)-\{w,w_1,\ldots,w_h\}$ such that $|E(z)\cap F_e|=0$. Otherwise, $|E(w)\cap F_e|=h-1$ for any vertex $w\in V-(F_1\cup F_2)$. Suppose $ww_i\in E(w)\cap F_e$ for $1\leq i\leq h-1$ and $W=\{w_i\mid 1\leq i\leq h-1\}$. If $W\cap(F_1\cup F_2)\neq\emptyset$, then $t-h+1=|F_1\cup F_2|\geq|F_1\cup F_2-W|+|W\cap(F_1\cup F_2)|\geq|N_{G-F_e}(w)|+1\geq\delta(G)-(h-1)+1>t-h+1$, a contradiction. Thus, $W\subseteq V-(F_1\cup F_2)$ and $|V-(F_1\cup F_2)-W-\{w\}|\geq(h+3)-(h-1)-1=3$. Suppose $z_i\in V-(F_1\cup F_2)-W-\{w\}$ for $1\leq i\leq3$. Then $|(\cup_{i=1}^3E(z_i))\cap(F_e-E(w))|\leq1$. Hence, $|E(z_j)\cap F_e|=0$ for some $j\in\{1,2,3\}$. Since $0=|E(z_j)\cap F_e|\geq h-1$, we have $0\leq h\leq1$. 

If $h=0$, then $|F_1\cap F_2|=t-1$ and $|V-(F_1\cup F_2)|=|V|-|F_1\cup F_2|\geq2t+3-(t+1)=t+2$. Pick a vertex $x'\in V-(F_1\cup F_2)$. Note that $N_{G}(x')\subseteq F_1\cup F_2$ by Claim 1. Thus, $t+1=|F_1\cup F_2|\geq |N_{G}(x')|\geq\delta(G)\geq t$.
If $\delta(G)=t$, then $|F_1\cup F_2|=\delta(G)+1$ and $|V-(F_1\cup F_2)|\geq\delta(G)+2$.
Hence, $G$ is isomorphic to $\Gamma_5(\delta(G), l)$ for some graph $H_{\delta(G)+1}$ and some integer $l$, a contradiction.
If $\delta(G)=t+1$, then $|F_1\cup F_2|=\delta(G)$ and $|V-(F_1\cup F_2)|\geq\delta(G)+1$. Hence, $G$ is isomorphic to $\Gamma_1(\delta(G), l)$ for some graph $H_{\delta(G)}$ and some integer $l$, a contradiction.

If $h=1$, then $|F_1\cap F_2|=t-2$, $|F_1\cup F_2|=t\leq\delta(G)$ and $|V-(F_1\cup F_2)|=|V|-|F_1\cup F_2|\geq2(t-1)+3-t=t+1\geq3$. Thus, there exists a vertex $z'\in V-(F_1\cup F_2)$ such that $F_e\cap E(z')=\emptyset$. Note that $N_{G-F_e}(z')=N_{G}(z')\subseteq F_1\cup F_2$ by Claim 1. We have $|F_1\cup F_2|\geq\delta(G)$. Therefore, $t=\delta(G)$. Without loss of generality, we assume that $F_e=\{u_1u_2\}$.

If $F_e\cap E(G[F_1\cup F_2])\neq\emptyset$, then $|N_{G-F_e}(u)|=|N_G(u)|\geq\delta(G)=|F_1\cup F_2|$ for any vertex $u\in V-(F_1\cup F_2)$. Since $N_{G-F_e}(u)\subseteq F_1\cup F_2$, we have $|N_{G}(u)|=\delta(G)$ for any vertex $u\in V-(F_1\cup F_2)$. By Claim 1, $G$ is isomorphic to $\Gamma_1(\delta(G), l)$ for some graph $H_{\delta(G)}$ and some integer $l$, a contradiction.

If $F_e\cap E[F_1\cup F_2, V-(F_1\cup F_2)]\neq\emptyset$, then we can suppose that $u_1\in V-(F_1\cup F_2)$ and $u_2\in F_1\cup F_2$. Note that $N_{G-F_e}(u_1)\subseteq F_1\cup F_2-\{u_2\}$. Then $\delta(G)-1\leq |E(u_1)|-|F_e|=|N_{G-F_e}(u_1)|\leq |F_1\cup F_2-\{u_2\}|=\delta(G)-1$. Therefore, $|N_G(u_1)|=\delta(G)$ and
$G$ is isomorphic to $\Gamma_1(\delta(G), l)$ for some graph $H_{\delta(G)}$ and some integer $l$, a contradiction.

If $F_e\cap E(G[V-(F_1\cup F_2)])\neq\emptyset$, then $u_1, u_2\in V-(F_1\cup F_2)$. Note that $N_{G-F_e}(u_i)\subseteq F_1\cup F_2$ for $i\in\{1,2\}$. Then $\delta(G)-1\leq |E(u_i)|-|F_e|=|N_{G-F_e}(u_i)|\leq |F_1\cup F_2|=\delta(G)$.
Therefore, $|N_G(u_i)|\in\{\delta(G), \delta(G)+1\}$ for $i\in\{1,2\}$ and $G$ is isomorphic to $\Gamma_4(\delta(G), l)$ for some graph $H_{\delta(G)}$ and some integer $l$ (see Figure \ref{11}(d)), a contradiction.

As mentioned above, we complete the proof of Theorem \ref{LL3}. $\qed$

Let $G$ be a $t$-regular $t$-connected graph and $h=0$. By Lemma \ref{Pro} (1) and Theorem \ref{LL3}, we have the following result.
\begin{corollary}[\cite{ChaC}]\label{cor5}
Let $G=(V, E)$ be a $t$-regular $t$-connected network with $N$ nodes and $t>2$. $G$ is $t$-diagnosable under the MM$^*$ model if $N\geq 2t+3$.
\end{corollary}

Note that a maximally connected graph $G$ is $\delta(G)$-connected. By Lemma \ref{LL1} and Theorem \ref{LL3}, we obtain the following result.
\begin{thm}\label{main2}
Let $G=(V, E)$ be a maximally connected graph with $|V|\geq 2(\delta(G)-h)+3$ and $\delta(G)\geq3$. If $G\notin \mathcal{F}(\delta(G))$, then $t_h^e(G)=\delta(G)-h$ under the MM$^*$ model for $0\leq h \leq \lfloor\dfrac{\delta(G)-1}{2}\rfloor$.
\end{thm}

Note that a $k$-regular $k$-connected graph is a maximally connected graph. By Lemma \ref{Pro} (1) and Theorem \ref{main2}, we immediately obtain the following result.
\begin{corollary}\label{cor2}
Let $G=(V, E)$ be a $k$-regular $k$-connected graph with $|V|\geq 2(k-h)+3$ and $k\geq3$. Then $t_h^e(G)=k-h$ under the MM$^*$ model for $0\leq h \leq \lfloor\dfrac{k-1}{2}\rfloor$.
\end{corollary}

By Lemma \ref{Pro} (2), Remark \ref{r1} and Theorem \ref{main2}, we can also obtain the following result.
\begin{corollary}\label{cor3}
Let $G=(V, E)$ be a maximally connected graph with $|V|\geq 2(\delta(G)-h)+3$. If $\delta(G)\geq3$ and $C(G)\leq \delta(G)-2$ (or $\delta(G)\geq4$ and $C(G)\leq \delta(G)-1$), then $t_h^e(G)=\delta(G)-h$ under the MM$^*$ model for $0\leq h \leq \lfloor\dfrac{\delta(G)-1}{2}\rfloor$.
\end{corollary}

\section{Conclusions}\label{6}
In this paper, we obtain the lower bound of the $h$-edge tolerable diagnosability of a $t$-connected graph and establish the $h$-edge tolerable diagnosability of a maximally connected graph under the PMC model and the MM$^*$ model, which extends some results in \cite{ChaC,Hak,Lian}. By our main results, the $h$-edge tolerable diagnosabilities of many well-known irregular networks can be determined under the PMC model and the MM$^*$ model.
\section*{Acknowledgement}
Y. Wei's research is supported by the Natural Science Foundation of Shanxi Province (No. 201901D211106). W. Yang's research is supported by the National Natural Science Foundation of China (No. 11671296).

\end{document}